\documentclass[12pt]{article}
\usepackage{amssymb,amsmath,amsfonts,amsthm}
\newcounter{parag}

\newtheorem{lem}{Lemma}
\newtheorem{theorem}{Theorem}

\newtheorem{question}{Question}

\begin{document}

\begin{center}
{\bf \Large On characterisation of a finite group by the set of conjugacy class sizes}

\medskip

{\bf Ilya Gorshkov}
\footnote{The work is supported by  Russian Science Foundation (project
19-71-10016).}
\end{center}
{\it Abstract: Let $G$ be a finite group and $N(G)$ be the set of its conjugacy class sizes. 
In the 1980's Thompson conjectured that the equality $N(G)=N(S)$, where $Z(G)=1$ and $S$ is simple, implies the isomorphism $G\simeq S$. In a series of papers of different authors Thompson's conjecture was proved. In this paper, we show that in some cases it is possible to omit the conditions $Z(G)=1$ and $S$ is simple and prove a more general result.

\smallskip

 \smallskip
}

\section*{Introduction}

Consider a finite group $G$. For $g\in G$, denote by $g^G$ the conjugacy class of $G$ containing $g$, and by $|g^G|$ the size of $g^G$. The centralizer of $g$ in $G$ is denoted by $C_G(g)$. Denote by  $N(G)=\{|g^G|\ |~ g\in G \}$. In
1987 John Thompson posted the following conjecture concerning $N(G)$.
\medskip

\textbf{Thompson's Conjecture (see \cite{Kour}, Question 12.38)}. {\it If $L$ is a finite
simple non-abelian group, $G$ is a finite group with trivial
center, and $N(G)=N(L)$, then $G\simeq L$.}

\medskip

In a series of papers \cite{Ah16}, \cite{Ah17}, \cite{Gor2}, Thompson's conjecture was studied for different groups, and it was finally proved in the paper \cite{Gor} in 2019. We say that a group $L$ is recognizable by the set of conjugacy class sizes among finite groups with trivial center (briefly recognizable) if the equality $N(L)=N(G)$, where $G$ is a finite group with trivial center, implies the isomorphism $L\simeq G$. Since $N(L)=N(L\times A)$ for any abelian group $A$, the condition $Z(L)=1$ is essential but not necessary. Given a finite group $G$, consider pairs of groups $(\Gamma, Z)$, such that $Z\leq Z(\Gamma)\cap \Gamma'$ and $\Gamma/Z\simeq G$. The largest by order possible second component of a pair $(\Gamma,Z)$ associated with a given group $G$ is called the Schur multiplier of $G$, and it is denoted by $M(G)$. The following problem arises.

\begin{question}
Let $L$ be a non-abelian simple group and $H(L)=M(L).L$, where $M(L)$ is a Schur multiplier of $L$. Let $G$ be a group with the property $N(G)=N(H(L))$. It is true that $G\simeq H(L)\times A$, for abelian group $A$?
\end{question}

In this article we make the first step towards answering this question. We give a positive answer to this question in the case $L\simeq Alt_5$.

\begin{theorem}
If $G$ is a such that $N(G)=N(SL_2(5))$, then $G\simeq SL_2(5)\times A$, where $A$ is an abelian group.
\end{theorem}

It is easy to show that $Sym_3$ is recognizable. Thus, the condition of solvability is not a necessary condition for the recognizability. As an example of a non-recognizable group, we can take a Frobenius group of order $18$. There exist two non-isomorphic Frobenius groups of order $18$ with the same sets of conjugacy class sizes. Navarro \cite{Navar} found two finite groups $G$ and $H$ with trivial center such that $N(G)=N(H)$, where $G$ is solvable and $H$ is non-solvable. However, unrecognizable groups with a trivial center are very rare. In particular, all the known unrecognizable groups have non-trivial solvable radicals. We denote the direct product of $n$ copies of a group $G$ by $G^n$.
The following question generalises Thompson's conjecture.

\begin{question}
Let $S$ be a non-abelian simple group. Is it true that for every $n\in \mathbb{N}$ the group $S^n$ is recognizable?
\end{question}

In this paper we prove the following theorem.

\begin{theorem}
If $G$ is a group such that $N(G)=N(Alt_5\times Alt_5)$ and $Z(G)=1$, then $G\simeq Alt_5\times Alt_5$.
\end{theorem}

\section{Notations and preliminary results}
\begin{lem}[{\rm \cite[Lemma 1.4]{GorA2}}]\label{factorKh}
Let $G$ be a finite group, $K\unlhd G$ and $\overline{G}= G/K$. Take $x\in G$ and $\overline{x}=xK\in G/K$.
Then the following conditions hold

(i) $|x^K|$ and $|\overline{x}^{\overline{G}}|$ divide $|x^G|$.

(ii) If $L$ and $M$ are consequent members of a composition series of $G$, $L<M$, $S=M/L$, $x\in M$  and
$\widetilde{x}=xL$ is an image of $x$, then $|\widetilde{x}^S|$ divides $|x^G|$.

(iii) If $y\in G, xy=yx$, and $(|x|,|y|)=1$, then $C_G(xy)=C_G(x)\cap C_G(y)$.

(iv) If $(|x|, |K|) = 1$, then $C_{\overline{G}}(\overline{x}) = C_G(x)K/K$.

(v) $\overline{C_G(x)}\leq C_{\overline{G}}(\overline{x})$.
\end{lem}

\begin{lem}\label{NSL25}
$N(SL_2(5))=\{12,20,30\}$.
\end{lem}
\begin{proof}
The proof of this lemma is simple exercise.
\end{proof}

\begin{lem}\label{NA5A5}
If $\alpha\in N(Alt_5\times Alt_5)\cup\{1\}$, then $\alpha=a\cdot b$ where $a,b\in\{1, 12, 15, 20\}$.
\end{lem}
\begin{proof}
The proof of this lemma is simple exercise.
\end{proof}


\begin{lem}[{\rm\cite[Corollary 1]{Camin}}]\label{divisors}
If for some prime $p$ there is in $G$ no element whose index is divisible by $p$, then either $G$ has order prime to $p$, or $G=H\times K$, where $H$ has order prime to $p$ and $K$ is an abelian $p$-group
\end{lem}

\begin{lem}[{\rm \cite[Theorem 5.2.3]{Gore}}]\label{Gore5hzOcomutantePstavtomor}
Let $A$ be a $\pi(G)'$-group of automorphisms of
an abelian group $G$. Then $G=C_G(A)\times[G,A]$.
\end{lem}

\section{Proof of Theorem 1}

Let $G$ be a group such that $N(G)=N(SL_2(5))$. Let $G=T\times A$, where $A$ is an abelian group, and $T$ does not include abelian direct factors. From Lemma \ref{NSL25} it follows that $N(G)=\{12,20,30\}$. From this fact and Lemma \ref{divisors} it follows that $\pi(T)=\{2,3,5\}$. Let us prove that $T\simeq SL_2(5)$.

We will show that $G$ contains a $5$-element $x$ such that $|x^G|=12$. Let $g\in G$ be an element of minimal order with property $|g^G|=12$. Since $12$ is a minimal by divisibility in $N(G)$ and $|(g^a)^G|$ divides $|g^G|$ for any $a\in \mathbb{N}$, we have $|(g^a)^G|\in\{1,12\}$. From the minimality of $g$ it follows that $|g|=p^b$ for some $b\in \mathbb{N}$, where $p\in\{2,3,5\}$. Assume that $p\neq 5$. We have that $C_G(g)$ includes a Sylow $5$-subgroup $S_5$ of $G$. Take $h\in S_5$. We have $C_G(gh)=C_G(g)\cap C_G(h)$, in particular, $|g^G|$ and $|h^G|$ divide $|(gh)^G|$. Since $12$ is the maximal and the minimal by divisibility element of the set $N(G)$, we conclude that $|(hg)^G|= 12$. Therefore $C_G(h)\geq C_G(g)$, in particular, $C_G(h)\geq S_5$. Therefore we have $S_5$ is abelian. Since there exists an element $y\in G$ such that $|y^G|_5>1$, we see that $Z(T)$ does not include $S_5$. Therefore there exists $5$ element $x$-such that $|x^G|>1$. We know that $S_5$ is abelian. Hence $|x^G|_5=1$ and $|x^G|=12$.

Similar we can show that there exists a $3$-element $y\in G$ such that $|y^G|=20$. Since $N(T)$ does not contain a number divisible by $12$ and $20$, we have $x'\in x^G$ and $y'\in y^G $ do not commute.

Suppose that $T$ is solvable. Let $H$ be a Hall $\{3,5\}$-subgroup of $T$ such that $C_G(x)\cap H$ is a Hall $\{3,5\}$-subgroup of $C_G(x)$. Put $y'\in y^G\cap H$ and $C_G(y')\cap H$ is a Hall $\{3,5\}$-subgroup of $C_G(y')$. Therefore, we have $|x^H|=3$, $|(y')^H|=5$ and $xy'\neq y'x$. Also we can assume that $x$ and $y'$ have minimal orders among the elements with this properties. Let $T\lhd H$ be a maximal subgroup among subgroups which does not contain $x$ and $y'$. Put $\overline{H}=H/T$, $\overline{x}\in\overline{H}$ is the image of $x$, $\overline{y'}\in\overline{H}$ is the image of $y'$, $R$ is the minimal normal subgroup of $\overline{H}$. Minimality of $x$ and $y'$ implies that $\overline{x}$ and $\overline{y'}$ do not commute, in particular, $|\overline{x}^{\overline{H}}|=3$ and $|\overline{y'}^{\overline{H}}|=5$. From definition of $T$ it follows that $\overline{x}\in R$ or $\overline{y'}\in R$. Suppose that $\overline{x}\in R$. Therefore $\overline{y'}$ acts non-trivially on $R$. From Lemma \ref{Gore5hzOcomutantePstavtomor} it follows that $R= C_{\overline{H}}(\overline{y'})\times [\overline{H},\overline{y'}]$. The element $\overline{y'}$ acts freely on $[\overline{H},\overline{y'}]$. Therefore $|[\overline{H},\overline{y'}]|\geq 5^2$. From Lemma \ref{factorKh} it follows that $|[\overline{H},\overline{y'}]|$ divides $|\overline{y'}^{\overline{H}}|$ and we have the contradiction. Similar we can show that $R$ does not contain $\overline{y'}$. Therefore $T$ is non-solvable.

From the description of simple groups with  the order divisible by primes which are less then or equal $5$ see \cite{Zav} it follows that $T$ includes a composition factor isomorphic $Alt_5$. Let $R\lhd T$ be a maximal subgroup with the property $F\leq T/R$ where $F\simeq Alt_5$. Using Lemma \ref{factorKh} we can show that $T/R=F$.

Let $H\lhd T$ be minimal subgroup. Therefore $H$ is a $p$-subgroup. Let $g\in T$ be a $p'$-element such that its image to $T/R$ is not trivial. We can show that $g$ acts trivially on $H$. Therefore $\langle g^T\rangle$ acts trivially on $H$. Using this reasoning, we can show that $R\leq Z(T)$. Since $T$ does not include direct factors, we see that $R=M(T)$. Hence, $T=M(F).F\simeq SL_2(5)$. Theorem 1 is proved.

\section{Proof Theorem 2}

Let $G$ be a group such that $N(G)=N(Alt_5\times Alt_5)$ and $Z(G)=1$.

\begin{lem}\label{l12}
There exists a $5$-element $a\in G$ such that $|a^G|=12$.
\end{lem}
\begin{proof}
 Suppose that $G$ does not contain a $5$-element whose order of the conjugacy class is equal to $12$. Let $P_5\in Syl_5(G)$ and $a\in Z(P_5)$. Therefore $|a^G|=12\cdot12$. Take $b\in G$ an element of minimal order such that $|b^G|=12$. Hence $|b|=p$, where $p\in\{2, 3\}$. To simplify the notation, we can assume that $p=2$. Since $C_G(b)$ includes a Sylow $5$-subgroup, we can assume that $a\in C_G(b)$.
Put $P_3\in Syl_3(C_G(b))$ and $c \in Z(P_3)$. Therefore $|(bc)^G| \in\{12, 12\cdot20\}$. Assume that $|(bc)^G|=12$. Consequently $C_G(b)=C_G(c)$. This statement implies that $|x^G|=20$ for each element $x$ from the centre of a Sylow $3$-subgroup of $G$. In particular, $C_G(c)$ contains a $3$-element $f$ such that $|f^G|=20$. Therefore, $f\in C_G(b)$ and $|(bf)^G|>12$. Thus, we can assume that $|(bc)^G| = 12\cdot 20$. Since $|G|_5 \geq 25$, $C_G(bc)$ contains a $5$-element $h$.
We have $|h^G| \neq 12$. Therefore, $|(bch)^G|=12\cdot20=|h^G|$, which implies $C_G(bc)=C_G(h)$. Since $a$ centralizes some Sylow $5$-subgroup of $G$, it can be assumed that $a\in C_G(h)$. Thus $a\in C_G(bc)$ and therefore $|(bch)^G|$ is a multiple of $12\cdot12$ and $12\cdot 20$, what gives a contradiction.

\end{proof}

\begin{lem}\label{l125}
If $|x^G|=12$ then $x$ is a $5$-element.
\end{lem}
\begin{proof}
From Lemma \ref{l12} it follows that there exists a $5$-element $a$ such that $|a^G|=12$. Assume that there exists an element $b$ such that $|b^G|=12$ and $\pi(b)\neq \{5\}$. From the minimality $|b^G|$ we can think that $b$ is an element of prime power order. Assume that $|b|=2$. Since $C_G(b)$ includes a Sylow $5$-subgroup of $G$, we can assume that $|(ab)^G|_5=1$. Therefore $|(ab)^G|\in\{12,12\cdot12\}$.

Assume that $|(ab)^G|=12$. Therefore $C_G(a)=C_G(b)$. Let $d\in G$ be an element of minimal order such that $|d^G|=15\cdot 15$. The subgroup $C_G(d)$ there includes a Sylow $2$-subgroup of $G$. Therefore we can assume that $b\in C_G(d)$. If $2\not \in \pi(d)$ then $C_G(bd)=C_G(b)\cap C_G(d)$, in particular $|(bd)^G|$ is a multiple of $12\cdot5\cdot5$; a contradiction. Since $C_G(a)=C_G(b)$, we can assume that $\{2,5\}\subseteq\pi(d)$. Let $d=d_2d_3d_5$, where $d_2=d^{|d|_{2'}}, d_3=d^{|d|_{3'}}, d_5=d^{|d|_{5'}}$. Assume that $d_3\neq 1$. Since $|d|$ is minimal, we see that $|(d_2d_3)^G|=|(d_2d_5)^G|=|(d_5d_3)^G|=15$, in particular, $C_G(d_2)=C_G(d_3)=C_G(d_5)$. Hence $C_G(d)=C_G(d_2)\cap C_G(d_3)\cap C_G(d_5)=C_G(d_2)$ and therefore $|d^G|=15$; a contradiction. Therefore $d=d_2d_5$. Similar we can show that $C_G(d)$ does not contain an element of order $3$. Hence $|G|_3=|d^G|_3=9$. Thus the Sylow $3$-subgroups of $G$ is abelian. We have $|G|_5=|d^G|_5|C_G(d)|_5\geq25|d|_5>25$.

Let $r\in G$ be an element of minimal order such that $|r^G|=20\cdot 20$. If $r$ is a $5$-element, then $|(rb)^G|$ ia a multiple of $12\cdot 5\cdot 5$; a contradiction. Therefore $\pi(r)\neq\{5\}$. From the fact that $|G|_5>25$, it follows that $C_G(r)$ contains an element $t$ of order $5$. Since $|r^G|$ is maximal by divisibility, $|t^G|$ divides $|r^G|$. Therefore $|t^G|=20$. Put $h\in Z(P_3(C_G(b))$, where $P_3(C_G(b))\in Syl_3(C_G(b))$. Therefore $|(bh)^G|=12\cdot 20$, in particular $|h^G|=20$. We have $|(bh)^G|_5=|h^G|_5$. Hence $C_G(bh)$ includes a Sylow $5$-subgroup of $C_G(h)$. The element $t$ centralize some Sylow $3$-subgroup of $G$. We can think that $h\in C_G(t)$. Since $C_G(b)$ there includes some Sylow $5$-subgroup of $C_G(h)$ it follow that $C_G(b)$ there contains $t$. Therefore, $|(ht)^G|$ divide $|(bht)^G|=12\cdot 20$. Since a Sylow $3$-subgroup of $G$ is abelian it follow that $|(ht)^G|_3=1$. Therefore $|(ht)^G|=20$ and $C_G(h)=C_G(t)\cap C_G(h)=C_G(t)$. Let $f\in C_G(t)$ be such that $|f^G|=12\cdot12$. That is obvious that $C_G(f)$ does not contain an elements of order $3$; a contradiction. We prove that $G$ does not contain a $\{2,5\}$-element $x$ such that $|x^G|=12$ and $x^{|x|_{2'}}\neq1, x^{|x|_{5'}}\neq1$. In particular, $|(ab)^G|=12\cdot 12$.

Take $c\in Z(P_3(C_G(b)))$. We have $|(bc)^G|_3=|b^G|_3=3$. Therefore $|(bc)^G|\in\{12, 12\cdot20\}$. Assume that $|(bc)^G|=12$. Then $C_G(b)=C_G(c)$. In particular, $C_G(b)$ contains an element $x\in Z(P_3)$ where, $P_3\in Syl_3(G)$ is such that $P_3\cap C_G(b)\in Syl_3(C_G(b))$. We have $|x^G|_3=1$. Therefore $|(bx)^G|=12\cdot20$. We can assume that $|(bc)^G|=12\cdot20$. The group $C_G(bc)$ contains a maximal subgroup $H$ of $P_5$ for some $P_5\in Syl_5(C_G(b))$. Therefore $Z(P_5)\cap H>1$. Let $h\in Z(P_5)\cap H$. We have $C_G(h)\cap C_G(b)> P_5$. Hence, $|(hb)^G|_5=1$, in particular, $|(hb)^G|=12\cdot12$. Thus $|(hbc)^G|$ is a multiple of $12\cdot12$ and $12\cdot 20$; a contradiction.

We prove that $G$ does not contain a $2$-element $b$ such that $|b^G|=12$. Similar we can show that $G$ does not contain a $3$-element $y$ such that $|y^G|=12$. Since $|(g^{\alpha})^G|$ divides $|g^G|$ for any $\alpha\leq|g|$, we obtain assertion of the lemma.
\end{proof}

\begin{lem}\label{l152}
If $|a^G|=15$, then $a$ is a $2$-element.
\end{lem}
\begin{proof}
Similar to Lemma \ref{l125}.
\end{proof}
\begin{lem}\label{l203}
If $|a^G|=20$, then $a$ is a $3$-element.
\end{lem}
\begin{proof}
Similar to Lemma \ref{l125}.
\end{proof}
\begin{lem}\label{l1212}
If $|a^G|=12\cdot 12$, then $a$ is a $5$-element.
\end{lem}
\begin{proof}
Assume that there exists an element $b$ such that $|b^G|=12\cdot12$ and $\pi(b)\neq\{5\}$. From Lemma \ref{l12} it follows that $|(b^{\alpha})^G|=|b^G|$ for each $\alpha$ such that $|b|/\alpha$ is not a power of $5$. In particular, there exists $\alpha$ such that $|b^{\alpha}|=p$ where $p\in\{2,3\}$, and $|(b^{\alpha})^G|=12\cdot 12$. Let $p=2$. Denote by $x=b^{\alpha}$. Using Lemma \ref{l152} we can show that there exists an element $h$ of prime power order such that $|h^G|=15\cdot 15$. Since $C_G(h)$ contains some Sylow $2$-subgroup of $G$, we can assume that $x\in C_G(h)$. If $\pi(h)\neq\{2\}$, then $C_G(xh)=C_G(x)\cap C_G(h)$, in particular, $|(xh)^G|$ is a multiple of $|x^G|$ and $|h^G|$; it contradicts the fact that $|x^G|$ and $|h^G|$ are distinct and maximal by divisibility. Therefore $h$ is a $2$-element. Similar as before we can show that for some $t\in\{2,3,5\}$ there exists a $t$-element $g\in G$, such that $|g^G|=20\cdot20$. Assume that $t=2$. Since $C_G(h)$ includes a Sylow $2$-subgroup of $G$, we can assume that $h\in C_G(g)$. The group $C_G(g)$ includes some Sylow $3$-subgroup of $G$. Lemma \ref{l203} implies that $C_G(g)$ contains a $3$-element $c$ such that $|c^G|=20$. We have $|(cg)^G|=20\cdot 20=|g^G|$. Therefore $C_G(c)> C_G(g)$, in particular $h\in C_G(c)$. Thus $|(ch)^G|$ is a multiple of $|c^G|$ and $|h^G|$; it contradicts the fact that $|h^G|$ is maximal by divisibility in $N(G)$. Therefore $t\in\{3,5\}$. Since $|G|_2>|b^G|_2=|c^G|_2$, we see that $C_G(c)$ contains an element $h'$ of order $2$. From Lemma \ref{l203} it follows that $|(h')^G|=20\cdot20$. Similar as above accept the contradiction. Therefore $p\neq2$.
Similar we can show that $p\neq 3$.
\end{proof}
\begin{lem}\label{l2015}
If $|a^G|=20\cdot 20$, then $\pi(a)=\{3\}$.
If $|a^G|=15\cdot 15$, then $\pi(a)=\{2\}$.
\end{lem}
\begin{proof}
Similar to Lemma \ref{l1212}.
\end{proof}

\begin{lem}
$|G|=60^2$.
\end{lem}
\begin{proof}
Let $a\in G$ be such that $|a^G|=12\cdot12$. From Lemma \ref{l1212} it follows that $|a|$ is a power of $5$. Since $|a^G|$ is maximal, it follows that for each $5'$-element $x\in C_G(a)$ the order of $x^G$ divides $|a^G|$. Therefore $|x^G|\in\{12, 12\cdot12\}$. Lemmas \ref{l125} and \ref{l1212} implies that $C_G(a)$ is a $5$-group. Thus $|G|_{5'}=|a^G|_{5'}=9\cdot 16$. Let $b\in G$ be such that $|b^G|=15\cdot15$. Using Lemmas \ref{l152} and \ref{l2015} we can show that $|G|_5=|b^G|_5=25$. Since $\pi(G)=\{2,3,5\}$, it follows that $|G|=60^2$.
\end{proof}
\begin{lem}
$G\simeq Alt_5\times Alt_5$.
\end{lem}
\begin{proof}
Let $L$ be a minimal normal subgroup of $G$. Suppose that $L$ is a $2$-group. Assume that there exists an element $a\in L$ such that $|a^G|=15$. Then $|C_G(a)|=16\cdot3\cdot 5$. Since Sylow $3$- and $5$-subgroups of $G$ are abelian, for any $3$-element $x\in C_G(a)$ we have $|x^G|=20$, in particular, $|(ax)^G|=15\cdot20$. For each $5$-element $y\in C_G(a)$ we have $|y^G|=12$, in particular, $|(ax)^G|=15\cdot12$. Suppose that $C_G(a)$ includes a Hall $\{3,5\}$-subgroup $H$. Since $|H|=15$, we see $H=\langle b\rangle $ is a cyclic group. Hence, $|(ab)^G|$ is a multiple of $|(ab^3)^G|=15\cdot 12$ and $|(ab^5)^G|=15\cdot 20 $; a contradiction. Thus, $C_G(a)$ is not solvable. Let $X$ be a minimal non-solvable subgroup of $C_G(a)$. Suppose that $X$ acts non-trivial on $L$. Then the automorphism group of the group $L/\langle a \rangle$ is non-solvable. Therefore, $|L/\langle a \rangle| \geq 2^4 $. Since $L$ is abelian, we get that $|C_G(a)|_2 \geq 2^5 $; a contradiction. Thus $X$ acts trivially on $L$. Let $c\in G$ be such that $|c^G|=12\cdot 12$. From Lemma \ref{l1212}, it follows that $c$ is a $5$-element. We have $C_G (c)$ is a $5$-group. Therefore $c$ acts freely on $L$. Hence, $|L|-1$ is a multiple of $5$, which means that $|L|\geq16$. Hence $X/X\cap L$ is a $2'$-group; a contradiction with the fact that $X$ is non-solvable. We have $|a^G|\neq 15$ for each $a\in L$. Therefore, for each $a\in L$ we have $|a^G|=15 \cdot 15$. In particular, $C_G(a)$ is a $2$-group. Let $P$ be a Sylow $5$-subgroup of $G$. Therefore $P$ acts freely on $L$. So $|L|>25 $; It contradicts the fact that the order of the Sylow $2$-subgroup of $G$ is equal $16$.

Suppose that $L$ is a $3$-group. Therefore $|L|\leq 9$. The automorphism group of the group $L$ does not contain an element of order $5$. Therefore any element of order $5$ acts trivially on $L$. Let $a \in G$ be such that $a^G=12\cdot 12$. From Lemma \ref{l1212} it follows that $a$ is a $5$-element. We have $|C_G(a)|=25$ and $L<C_G(a)$; a contradiction.


Suppose that $L$ is a $5$-group. If $|L|=5$, then the Sylow $3$-subgroup of $G$ acts trivial on $L$. Therefore, the centralizer of each $3$-element contains an element of order $5$. By Lemma \ref{l2015} we have that $G$ contains a $3$-element $b$ such that $|b^G|=20\cdot20$, in particular, $C_G(b)$ is a $3$-group; a contradiction. Thus $|L|=25$. In this case we have $G/L$ is a $\{2,3\}$-group, in particular, $G$ is solvable. Let $x\in G$ be such that $|x^G|=15$. We have that the Hall $\{3, 5\}$-subgroup of $C_G(x)$ must be a Frobenius group with the kernel of order $5$ and the complement of order $3$; a contradiction.


Thus $L$ is a direct product of simple groups. Assume that $L$ is not a simple group. The minimal simple group has order $60$ and it is isomorphic to $Alt_5$. We get that $G=L$. In this case, $L$ is not a minimally normal subgroup.


We proved that every minimal normal subgroup of $G$ is a non-abelian simple group. Let $C$ be the socle of $G$. We have $G\leq Aut(C)$. Assume that $C$ is a simple group. From the description of simple groups, of orders which prime divisors from $\{2,3,5 \}$ see \cite{Zav}, we see that $C$ is isomorphic to one of the groups $Alt_5, Alt_6, U_4 (2)$. Note that in this case $|C|_5 = 5<|G|_5$ and $|Out(C)|_5=1$; a contradiction.


Thus $C$ is a product of several simple groups. Similar as above we conclude that $C\simeq Alt_5\times Alt_5\simeq G$.


\end{proof}

\bigskip

Ilya~B. Gorshkov

Sobolev Institute of Mathematics SB RAS

Novosibirsk, Russia

E-mail address: ilygor8@gmail.com

\end{document}